# Design of a Fractional Order Phase Shaper for Iso-damped Control of a PHWR under Step-back Condition

Suman Saha[1], Saptarshi Das[1], Ratna Ghosh[2], Bhaswati Goswami[2], R. Balasubramanian[3], A. K. Chandra[3], Shantanu Das[4], Amitava Gupta[1]

*Abstract*—Phase shaping using fractional order (FO) phase shapers has been proposed by many contemporary researchers as a means of producing systems with iso-damped closed loop response due to a stepped variation in input. Such systems, with the closed loop damping remaining invariant to gain changes can be used to produce dead-beat step response with only rise time varying with gain. This technique is used to achieve an active *step-back* in a Pressurized Heavy Water Reactor (PHWR) where it is desired to change the reactor power to a pre-determined value within a short interval keeping the power undershoot as low as possible. This paper puts forward an approach as an alternative for the present day practice of a passive step-back mechanism where the control rods are allowed to drop during a step-back action by gravity, with release of electromagnetic clutches. The reactor under a step-back condition is identified as a system using practical test data and a suitable Proportional plus Integral plus Derivative (PID) controller is designed for it. Then the combined plant is augmented with a phase shaper to achieve a dead-beat response in terms of power drop. The fact that the identified static gain of the system depends on the initial power level at which a step-back is initiated, makes this application particularly suited for using a FO phase shaper. In this paper, a model of a nuclear reactor is developed for a control rod drop scenario involving rapid power reduction in a 500MWe Canadian Deuterium Uranium (CANDU) reactor using AutoRegressive Exogenous (ARX) algorithm. The system identification and reduced order modeling are developed from practical test data. For closed loop active control of the identified reactor model, the fractional order phase shaper along with a PID controller is shown to perform better than the present Reactor Regulating System (RRS) due to its iso-damped nature.

1. S. Saha, S. Das and A. Gupta are with Department of Power Engineering, Jadavpur University, Kolkata-700098, India (e-mail: the_suman@yahoo.co.in).
2. R. Ghosh and B. Goswami are with Department of Instrumentation and Electronics Engineering, Jadavpur University, Kolkata-700098, India.
3. R. Balasubramanian and A. K. Chandra are with R & D, Electronic System, Nuclear Power Corporation of India Ltd., Mumbai-400085, India.
4. S. Das is with Reactor Control Division, Bhaba Atomic Research Centre, Mumbai-400085, India.

## I. INTRODUCTION

BULK power reduction in a nuclear reactor is done under load following operations [1] or under some abnormal operating conditions. Power reduction is mainly done with the help of control rods. When power generated by a reactor is to be reduced, the control rods are inserted to a specified level and gradually the set point of the demand power is also reduced. Control rod insertion for rapid power reduction of a reactor is known as reactor step-back. Modern PHWRs with passive safety features rely on gravity to drive the control rods which are held above the reactor by electromagnetic clutches. These are released during a step-back and also automatically, in case of a power failure. In this system, the clutches are de-energized to release the rods and then again energized to grip them once the rods have dropped by a designated amount. But this passive step-back mechanism in addition creates a power undershoot while also producing a very sluggish response which is not desired at all because the safety constraints do not permit excessive power undershoot/overshoot in the nuclear reactor as reported in [2].

The approach presented in this paper assumes an active, motor driven step-back mechanism in which the motor speed is varied by PID controller in conjunction with a FO phase shaper. For this closed loop system, the change in control rod position (due to step-back) is the input and the actual global reactor power is the controlled variable. In order to do this, the reactor must be identified using the dynamics of power variation during a step-back. However, the major drawback of an active step back is that it lacks the essential safety feature for total power failure of a nuclear power plant. But, this problem can be overcome by putting shut-off rods within the reactor for ensuring safety issues, as reported in [3].

From the point-kinetic governing equations [1], [2] of a nuclear reactor, it is evident that it is a highly nonlinear system. So, the modeling of the reactor as a linear system for controller design is not so easy as the transfer function of the identified system changes depending on the operating conditions i.e. the initial power of the reactor and level of control rod insertion [4], [5]. Attempts are made to design PID controller for nonlinear reactor model in [1], [4], [5] considering different linearized transfer function models of the same nonlinear plant under different operating conditions. The methodology proposed in [1] designs different controllers for change in operating point and proposes the switching of controllers to cope with the nonlinear nature of the reactor. This is quite obvious that designing different controllers for different linearized models of a nonlinear system is not a feasible solution as far as a wide variation in operating point (due to the change in initial power, rod position or both) are concerned. In this paper, a robust FO phase shaper along with a Linear Quadratic Regulator (LQR) [6] tuned PID controller is proposed. This is expected to produce satisfactory closed loop response for a wide range of variation in operating point, with respect to that designed in [4], [5] due to enhanced parametric robustness.

For this purpose, the nonlinear reactor model is reduced to two standard First Order Plus Time Delay (FOPTD) and Second Order Plus Time Delay (SOPTD) models for the ease of controller synthesis. It is found that only the gain of the reduced order models get changed due to shift in operating point. Next, a LQR based PID tuning methodology [6] is applied to obtain a dead-beat power drop response under a step-back using the plant with minimum gain. The plant with the PID controller is then augmented with a FO phase shaper, to ensure iso-damping property for a considerable range of gain variation [7], so that such change in linearized plant transfer function can be easily handled. Thus, design of different controllers and their switching as per shift in operating point [1], can be avoided.

In the earlier work [8], an optimal linear state feedback regulator has been proposed for controlling the reactor power with change in reactivity. The present approach followed in this paper first identifies the nuclear reactor as FOPTD and SOPTD models and then uses a phase shaper in conjunction with a LQR tuned PID controller in the RRS. The reactor with its power regulator in closed loop is taken as the system and is controlled by a master controller and a phase shaper in closed loop with the master controller set point acting as the local set point. The phase shaper is designed using a flat phase criterion based on Bode's integral [9], [10], [11]. The methodology uses a constrained optimization technique that flattens the phase curve for a maximum frequency spread around gain cross-over frequency $\omega_{gc}$ of the system represented by the plant and its PID controller. The iso-damped response of the corresponding closed-loop system due to a sudden rod drop under a step-back condition makes it particularly useful for the use in designing an active step-back system for a 500MWe electric CANDU type PHWR, where the system gain has been found to be dependent on initial power. Thus an active step-back system realized with the methodology presented in the paper can be used seamlessly over a considerable range of reactor operation with the power undershoot constant under a step-back condition (insensitive to change in system gain). The need for iso-damped response [7], [12], [13] of the reactor is not only to handle the changing transfer function of the nonlinear plant but also to allow an increase in the gain of the phase shaper to get a considerably faster control action than the gravity for control rod insertion through the viscous medium, while the power undershoot remains the same.

The methodology presented in the paper proposes the use of a FO differ-integrator for phase shaping and automatically establishes the integer approximated FO phase shaper as a rational transfer function, which is effective for a frequency spread around a specified frequency viz. the gain crossover frequency ($\omega_{gc}$). The phase shaper designed by this methodology produces the widest flat-phase region around $\omega_{gc}$ of the system comprising the plant and its controller, with the phase margin fixed above a specified value. Thus, the resultant closed-loop system exhibits iso-damped step response, with constant overshoot for a variation of system gain within a range. The methodology is demonstrated using the reduced order models of the reactor. The results have been repeated using higher order Carlson's approximation of FO phase shaper [7], [14]-[16] with little difference in simulated results, implying that a simple phase shaper approximated by first order transfer functions is sufficient.

The rest of the paper is organized as follows. Section II presents the methodology for reactor model identification using ARX model and reduced order modelling under different operating conditions. Concept of phase shaping using Bode's integral is presented in section III and the phase shaper design for the identified reactor models are presented in section IV. Next, simulation results for an active step-back for the FOPTD and SOPTD models are presented in section V, followed by the conclusion in section VI and the references.

## II. REACTOR MODEL DEVELOPMENT

### A. Overview of the Present Reactor Control Mechanism

In a typical 500MWe CANDU type PHWR, the reactor power is controlled by the reactor regulating system (RRS), consisting of mainly three components viz. control rods (CR), adjustor rods (AR) and zone control components (ZCC). The ARs are provided for fast startup of the reactor. The CRs are provided for coarse control and ZCCs for fine control of the power level. Normally CRs are kept fully out, ARs are kept fully in while the ZCCs are partially filled with water. The CRs are made of neutron absorber materials like Boron or Graphite etc. They are used for bulk power reduction in a reactor when it is required to suddenly reduce the power from a certain operating level. The ARs are generally kept within the reactor under normal operating conditions, which is an extra burden of the reactor, as it absorbs neutrons and reduces the reaction rate. When it is required to have a faster reaction rate in the reactor, ARs are taken out of the reactor gradually. This introduces some amount of positive reactivity in the reactor and thus power level is increased. The ZCCs are provided to have local control of power. This is achieved by varying the water level of those specific zones of the reactor.

The set-point of the RRS is specified at the desired level, which is called the demand power and the control loop error is corrected by the continuous measurement of the reactor power level or bulk power by the Self-Powered-Neutron-Detectors (SPND). So, the goal of the reactor control system is to minimize the effective power error (EPE), which is the sum of the difference between demand power and bulk power and the

## B. Identification using ARX model

This section presents the identification of the reactor along with its regulating system taking CR position as input and the bulk power as output using ARX model. For identification, the reactor is visualized as a system with CR position (fraction of total drop) as input and the actual power (in percentage of maximum power produced) as output. The identification is based on data obtained from operating Indian PHWRs provided by the Nuclear Power Corporation of India Ltd. (NPCIL). The simulation data is provided for 14 seconds with 0.1 second of sampling time. The data is shown in graphical form for the 30% rod drop case in Fig.1.

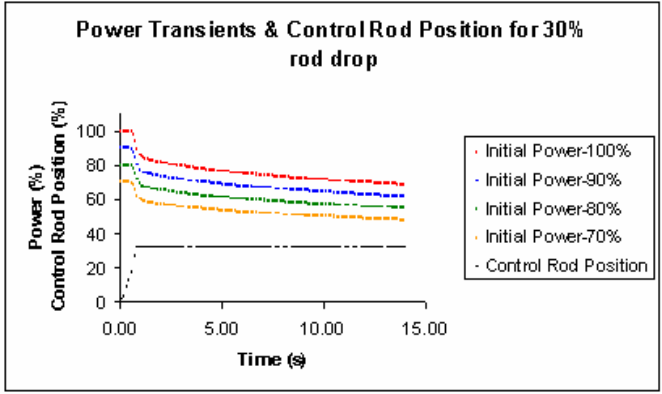

Fig. 1. Power Transients and Control rod position for rod drop upto 30%.

Now the basic identification using ARX models from a measured time domain data [17] is discussed. If it be assumed that at time $t$, the input and output of an unknown system is $u(t)$ and $y(t)$ respectively. Then the system can simply be described by the following linear difference equation (1)

$$y(t) + a_1 y(t-1) + \cdots + a_n y(t-n) = b_1 u(t-1) + \cdots + b_m u(t-m) \quad (1)$$

The above equation can be re-written in the following form if the values of input and output data at each time step are known:

$$y(t) = -a_1 y(t-1) - \cdots - a_n y(t-n) + b_1 u(t-1) + \cdots + b_m u(t-m) \quad (2)$$

The calculated value is thus

$$\hat{y}(t) = \varphi^T(t) \cdot \theta \quad (3)$$

where, $\theta = [a_1 \cdots a_n \; b_1 \cdots b_m]^T$ (4)

and $\varphi(t) = [-y(t-1) \cdots -y(t-n) \; u(t-1) \cdots u(t-m)]^T$ (5)

Using mean square error minimization over a time interval ($1 \le t \le N$), the system parameter vector $\theta$ can be approximated as

$$\hat{\theta} = \left[ \sum_{t=1}^{N} \varphi(t) \varphi^T(t) \right]^{-1} \sum_{t=1}^{N} \varphi(t) y(t) \quad (6)$$

Since the value of input and output at each instant i.e. $\varphi(t)$ vector is known, from (6) the coefficients of the discrete transfer function of the model i.e. $\hat{\theta}$ can be calculated. This is done using the *arx()* function of System Identification Toolbox of MATLAB and a discrete transfer function of an arbitrary order model has been estimated from measured input output data [18].

A 30% rod drop data is used to model a 500MWe CANDU reactor for different initial powers i.e. 100%, 90%, 80% and 70%. The ARX structure is used in this case to find out a model having minimum prediction error from the test data [17]. The identified transfer function models are found to be slightly different with each other since reactor dynamics is inherently non-linear in nature and its linearized model will differ depending on the shift in the operating point i.e. the initial power and level of control rod insertion [1]. The identified discrete transfer function models are then converted into continuous time models using Zero-Order Hold (ZOH) and minimum realization *minreal()* function of MATLAB's Control System Toolbox [19] (considering any possible pole zero cancellation) as follows:

$$G_1(s) = \frac{\begin{array}{l}(44.79s^5 + 8408s^4 + 7.687\times 10^4 s^3 + \\ +8.42\times 10^6 s^2 - 2.561\times 10^7 s + 1.336\times 10^7)\end{array}}{\begin{array}{l}(s^6 + 12.31s^5 + 1088s^4 + 6624s^3 + \\ +5.75\times 10^4 s^2 + 7.683\times 10^5 s + 6.946\times 10^4)\end{array}} \quad (7)$$

$$G_2(s) = \frac{\begin{array}{l}(-81.59s^5 + 8625s^4 - 2.028\times 10^4 s^3 + \\ +9.119\times 10^6 s^2 - 2.544\times 10^7 s + 1.682\times 10^7)\end{array}}{\begin{array}{l}(s^6 + 17.41s^5 + 1129s^4 + 9406s^3 + \\ +5.397\times 10^4 s^2 + 9.21\times 10^5 s + 9.474\times 10^4)\end{array}} \quad (8)$$

$$G_3(s) = \frac{\begin{array}{l}(22.75s^5 + 9232s^4 + 6.87\times 10^4 s^3 + \\ +7.943\times 10^6 s^2 - 2.047\times 10^7 s + 1.4\times 10^7)\end{array}}{\begin{array}{l}(s^6 + 14.49s^5 + 1101s^4 + 7680s^3 + \\ +5.278\times 10^4 s^2 + 8.547\times 10^5 s + 8.839\times 10^4)\end{array}} \quad (9)$$

$$G_4(s) = \frac{\begin{array}{l}(-61.92s^5 + 9106s^4 - 1.907\times 10^4 s^3 + \\ +7.272\times 10^6 s^2 - 2.017\times 10^7 s + 1.215\times 10^7)\end{array}}{\begin{array}{l}(s^6 + 15.31s^5 + 1105s^4 + 8861s^3 + \\ +5.144\times 10^4 s^2 + 9.169\times 10^5 s + 8.911\times 10^4)\end{array}} \quad (10)$$

The four transfer functions represent the reactor model for 30% rod drop case, having starting power level of 100%, 90%, 80% and 70% respectively. It is to be noted that all the transfer functions have one or more negative signs in the numerator, thus indicating non-minimum phase (NMP) time response of the system under all operating conditions.

## C. Sub-Optimal Model Reduction

The identified models, as discussed in the previous sub-section, are then reduced to simple FOPTD and SOPTD models for the ease of controller tuning by minimizing the difference between $H_2$-norm of the identified and reduced order model as presented in [20]. The reduced order models of the reactor at different starting power levels are described by (11)-(18).

The reduced order FOPTD models are as follows:

$$G_{100}^{I}(s) = \frac{384.6}{s+2}e^{-0.5s} \quad (11)$$

$$G_{90}^{I}(s) = \frac{355.1}{s+2}e^{-0.5s} \quad (12)$$

$$G_{80}^{I}(s) = \frac{316.7}{s+2}e^{-0.5s} \quad (13)$$

$$G_{70}^{I}(s) = \frac{272.6}{s+2}e^{-0.5s} \quad (14)$$

The subscripts with the transfer functions correspond to the initial power level at which the step-back is initiated and the superscripts correspond to the order of the reduced model. From equation (11)-(14), it is observed that the nonlinear nuclear reactor can be modelled as a generalized FOPTD transfer function of the structure $\frac{K}{sT+1}e^{-Ls}$ with varying system gain ($K$), while the time-constant ($T$) and the time-delay ($L$) remaining the same.

Similarly, the reduced order SOPTD models corresponding to equations (7)-(10) are identified as:

$$G_{100}^{II}(s) = \frac{192.3}{s^2+2s+1}e^{-0.5s} \quad (15)$$

$$G_{90}^{II}(s) = \frac{177.6}{s^2+2s+1}e^{-0.5s} \quad (16)$$

$$G_{80}^{II}(s) = \frac{158.4}{s^2+2s+1}e^{-0.5s} \quad (17)$$

$$G_{70}^{II}(s) = \frac{136.4}{s^2+2s+1}e^{-0.5s} \quad (18)$$

Proceeding in the same way, from equation (15)-(18), it is observed that the nonlinear nuclear reactor can also be modelled as a generalized SOPTD transfer function of the structure $\frac{K}{s^2+2\zeta\omega_n s+\omega_n^2}e^{-Ls}$ with system gain ($K$), while the damping-ratio ($\zeta$), undamped-natural frequency ($\omega_n$) and the time delay ($L$) remaining the same. The reduced order FOPTD and SOPTD models extract the dominant first order & second order modes in the system along with the inherent time delay which causes the NMP behaviour, while the linearly varying gain in these models are dependent on the operating point.

The advantage of the reduced order modelling of the reactor is that the nonlinear dynamic plant can easily be represented by a general class of linear varying gain system of FOPTD and SOPTD structure. Thus, with a robust tuning technique it is possible to maintain constant undershoot for variation of overall system gain (i.e. due to the nonlinearity of the plant and also a deliberate gain variation of the phase shaper to get faster time response) as presented in the next section.

### III. FORMULATION OF THE PHASE SHAPING METHODOLOGY

#### A. Philosophy of Phase Shaper Design

In this sub-section, the methodology of designing a phase shaper as a FO differ-integrator is presented. The detailed treatment of the design methodology can be found in [21]. For the phase shaper design, as stated in [21], the plant $G_{pl}(s)$ is assumed to be a FOPTD or a SOPTD system and the PID controller controlling the plant is assumed to be tuned by any standard method resulting in a stable closed loop system. Thus starting with the open loop system

$$G(s) = G_c(s) \times G_{pl}(s) \quad (19)$$

comprising the plant $G_{pl}(s)$ and its PID controller $G_c(s)$, a phase shaper is designed such that the resultant closed loop system exhibits iso-damped response to step changes in input over a range of gain variations. This is achieved by flattening the asymptotic phase curve of the system

$$G_{ol}(s) = G_{ph}(s) \times G_c(s) \times G_{pl}(s) \quad (20)$$

around its gain cross-over frequency, for the maximum possible frequency spread allowing maximum variation in system gain. This ensures a constant phase margin and hence gain-independent overshoot (iso-damped) for the time response of the system.

In (20), $G_{ph}(s)$ is the phase shaper realized using a FO differ-integrator

$$G_{ph}(s) = \frac{(1+as^q)}{s^q} \quad (21)$$

with $\frac{1}{a} \leq \omega_{gc}^q \quad (22)$

and $0 \leq q \leq 1 \quad (23)$

In (21) $\omega_{gc}$ represents the gain cross-over frequency of $G(s)$. The methodology put forward in [21] uses Bode's Integral to represent the phase of $G(s)$ around its gain cross-over frequency as

$$\left.\frac{d\angle G(j\omega)}{d\omega}\right|_{\omega=\omega_{gc}} = \frac{\phi_m - \pi}{\omega_{gc}} + \frac{2}{\pi\omega_{gc}}\ln|k_g| \quad (24)$$

$\phi_m$ is the phase margin of $G(s)$ and $k_g$ is its static gain. It is established in [9]-[11] that equation (24) is valid for both minimal and non-minimal phase systems alike. It is then attempted to flatten the phase around $\omega = \omega_{gc}$ using a phase shaper $G_{ph}(s)$, such that the condition

$$\frac{d}{d\omega}\angle G(j\omega) + \frac{d}{d\omega}\angle G_{ph}(j\omega) = 0 \quad (25)$$

is satisfied over a frequency band $\Delta\omega$, around $\omega_{gc}$.

Using (21) and (24), it follows from (25) that

$$\frac{\phi_m - \pi}{\omega_{gc}} + \frac{2}{\pi\omega_{gc}}\ln|k_g| + \frac{aq\omega^q \sin\frac{q\pi}{2}}{\omega(1+2a\omega^q \cos\frac{q\pi}{2}+a^2\omega^{2q})} = 0 \quad (26)$$

Thus, the addition of $G_{ph}(s)$ alters the phase of $G(s)$ and the net phase of $G_{ol}(s)$ at $\omega_{gc}$ can be expressed as

$$\phi\bigg|_{\omega=\omega_{gc}} = \phi_m - \pi - \frac{q\pi}{2} + \tan^{-1}\left(\frac{a\omega^q \sin\frac{q\pi}{2}}{1+a\omega^q \cos\frac{q\pi}{2}}\right) \quad (27)$$

It is seen from (24) that the phase of $G_{ol}(s)$ at $\omega_{gc}$ is less than the phase of $G(s)$ at the same frequency. Thus, the addition of the phase shaper flattens the phase curve of $G_{ol}(s)$ around $\omega_{gc}$, at the cost of reduction in phase margin. Thus if the minimum desired phase margin with the phase shaper be $\phi_{md}$, then it follows that the constraint

$$\phi_{md} - \phi_m + \frac{q\pi}{2} - \tan^{-1}\left(\frac{a\omega^q \sin\frac{q\pi}{2}}{1 + a\omega^q \cos\frac{q\pi}{2}}\right) \leq 0 \quad (28)$$

must be satisfied.

In [21], the problem of finding a phase shaper of the form represented by (21) that produces maximum flatness in terms of frequency spread is formulated as one of constrained optimization that finds $\{q, a\}$ maximizing the value of $|\omega - \omega_{gc}|$ and satisfying the constraints represented by equations (22), (23), (26) and (28) using MATLAB's Optimization Toolbox function *fmincon()* customized with an *Active set* algorithm [22]. The rationale behind use of *fmincon()* arises from the fact that this is an optimization problem with non-linear constraints and it has been shown that for such optimization problems *fmincon()* can be used effectively, as reported in [23]-[24], for example.

The methodology presented in [21] and used for the present work differs from the one proposed in [9]-[11] in the sense that the plant in this case is considered along with its PID controller, and thus it can be applied to any existing control loop without having to tune the controller again.

### B. Realization of FO Control Elements

A major issue in the use of FO controllers is its physical realization since fractional order elements are infinite dimensional linear filters. Vinagre *et al* have dealt with issues related to use of FOPID controllers for industrial applications in [25]. As reported in [25], a FO element can be realized by a lossy capacitor based micro-electronic approach [26], [27]; by analog circuit realization [28], [29] or by integer approximation using, for example, Carlson representation [7], [14]-[16] or the CRONE Toolbox [30]-[32]. Such integer approximations are valid for specific frequency ranges. For the present work, Carlson's representation is adopted which in the simplest case, approximating a FO differ-integrator (differentiator or integrator) as a rational transfer function with a first order numerator and a first order denominator. The present approach is independent of the realization method of FO elements.

## IV. DESIGN OF FO PHASE SHAPER FOR ISO-DAMPED RESPONSE

The reactor model for 70% initial power is found to have the minimum gain and it is then tuned with the LQR technique [6]. The significance of a PID controller tuned with a LQR technique can be understood by examining equation (27) which shows that the phase margin reduces with the introduction of the FO phase shaper. The methodology presented in [6] allows the LQR tuned PID controller, to be designed with a specified closed-loop damping and frequency with advantages of an optimal controller preserved. Thus if the PID controller is tuned to produce a slightly over-damped closed-loop response, and consequently a large phase margin, a FO phase shaper can be used to shape the phase curve at cross-over frequency to produce a dead-beat response while keeping a comfortable phase margin.

The LQR based tuning methodology as reported in [6] yields a PID controller for the FOPTD model as described in equation (14) as

$$G_c = 0.0059 + \frac{0.0019}{s} + 0.00082s \quad (29)$$

Using the methodology discussed in the previous section, the phase shaper for this controller is designed to be

$$G_{ph} = \frac{(1 + 1.3419 s^{0.6181})}{s^{0.6181}} \quad (30)$$

Fig. 2 shows the frequency response of the open-loop system which represents the reactor model (14) controlled by the existing RRS, the plant with controller (29) and the phase shaper (30). The flattening of the phase curve over a significant range of frequency is evident from a comparison of the phase curves. The increase in gain margin is also appreciable. Thus, if the system gain changes to values represented by equations (11) to (13), then the gain cross-over frequency of the plant with PID controller and phase shaper will shift to higher values, but the phase margin will remain constant. This implies that the corresponding closed-loop system damping will remain invariant to change in system gain, but the system rise time will decrease. Hence, if the design based on the slowest plant i.e. (14) with controller (29) and phase shaper (30) meets the timing constraints, then the closed loop system with the same controller and the phase shaper will meet the timing constraints with faster plants (11) to (13) also, while the damping remains constant.

It is clear that the phase shaper in conjunction with the PID controller makes phase curve of the overall system flat over a wide range of frequency, which ensures iso-damping property. This is evident from a comparison of the closed loop step-back responses (Fig. 3). As it can be observed from the Bode diagram (Fig. 2), the PID controller alone can not ensure flat phase around $\omega_{gc}$ and hence constant undershoot (Fig. 3).

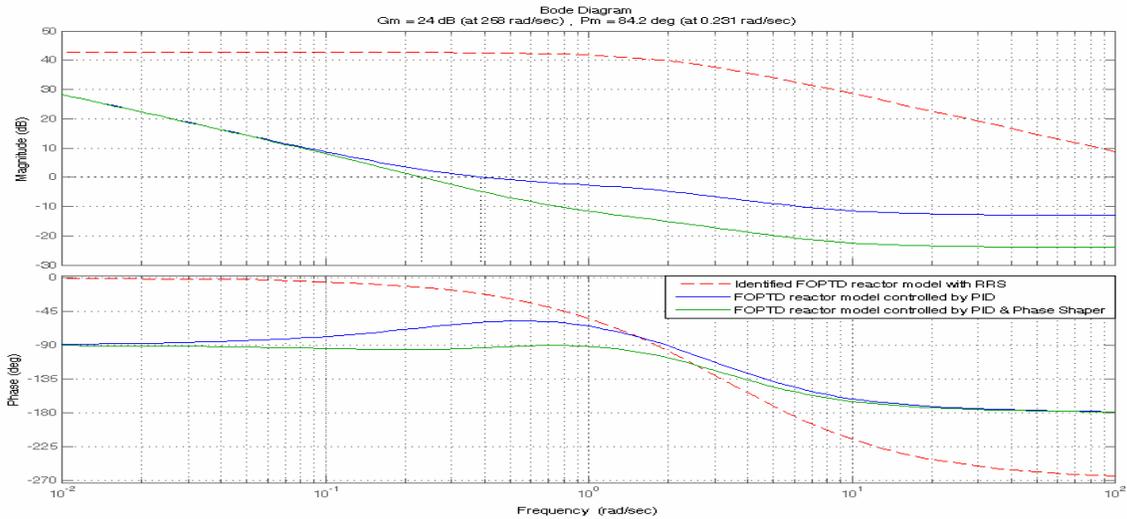

Fig. 2.  Phase flattening for the FOPTD reactor model (14) with PID controller (29) and phase shaper (30).

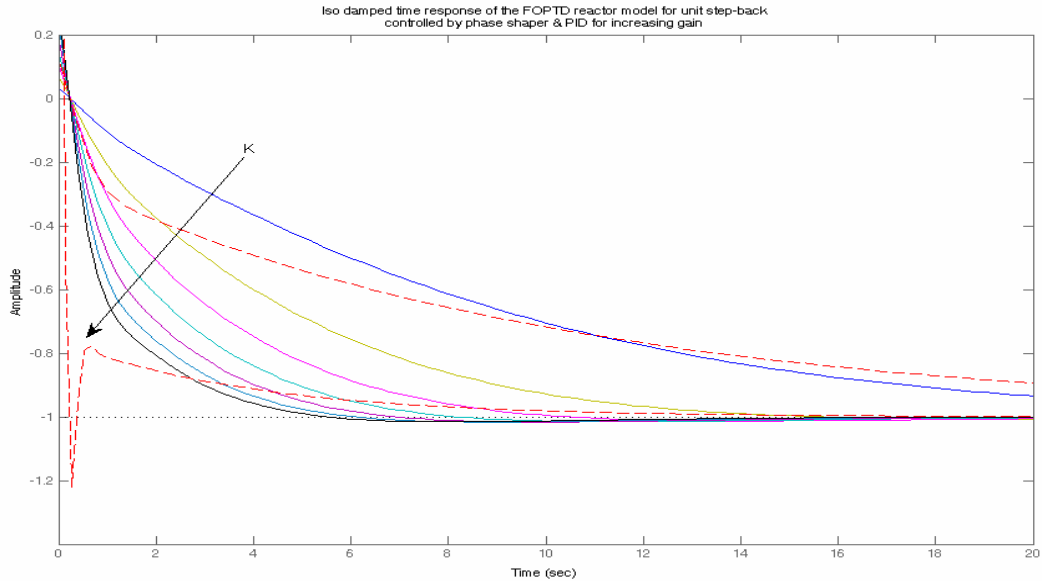

Fig. 3.  Iso damped time response of the FOPTD reactor model (14) for unit step-back signal; dotted lines represent an active step-back response with only PID controller (29) and continuous lines represent the response of the active rod drop with the PID along with the phase shaper (30) for 7 times increase in gain.

From Fig. 3, it is evident that with a simple PID, if the gain is increased to overcome the sluggish response, the closed loop response suffers from large undershoot. But using a phase shaper along with the PID, gives dead-beat response even with 7 times increase in the gain of the phase shaper, while giving much faster response and no undershoot. For a practical nonlinear system e.g. the nuclear reactor in our case, approximated as a Linear Time-Invariant (LTI) system, the system gain hardly exceeds 1.5-2 times. But the feasibility of 7 times increase in gain with the robust controller enables us to have a deliberate variation in the gain of the phase shaper to have much faster response without the occurrence of any undershoot.

Next, the methodology presented in section IIIA. is applied to obtain a PID controller and a corresponding FO phase shaper for the SOPTD approximated plant (18) as

$$G_c = 0.0039 + \frac{0.001}{s} + 0.002s \qquad (31)$$

and $G_{ph} = \dfrac{(1+5s^{0.75})}{s^{0.75}} \qquad (32)$

The flattening of the phase curve with the introduction of the controller (31) and phase shaper (32) as compared to the original plant represented by (18) and the plant with PID controller (31) is shown in Fig. 4. Clearly, the phase shaper in conjunction with the PID controller ensures iso-damping property for the SOPTD model also. This is evident in the closed loop step-down responses (Fig. 5).

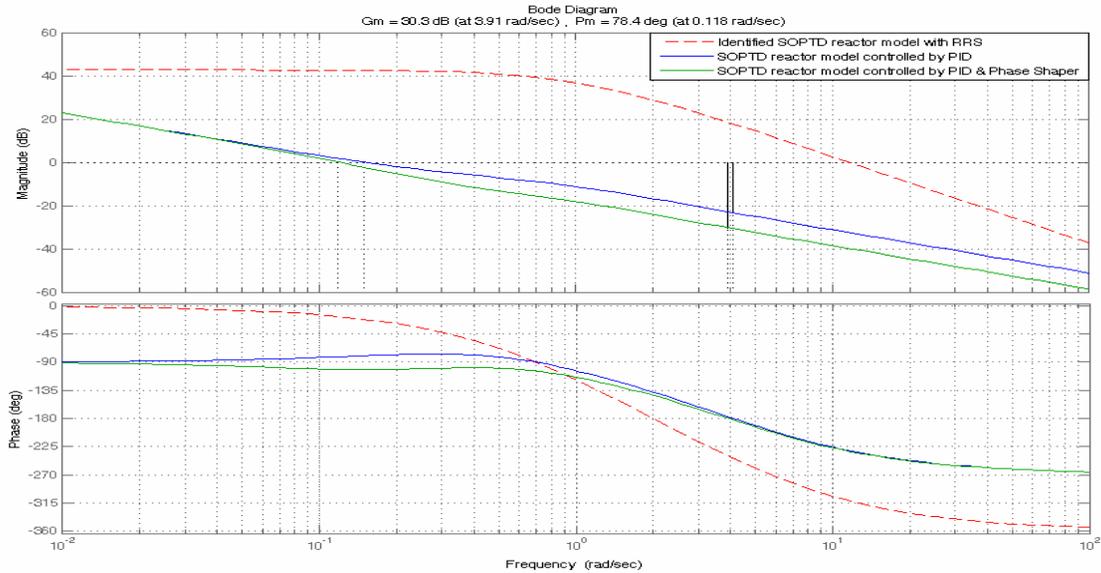

Fig. 4. Phase flattening for the SOPTD reactor model (18) with PID controller (31) and phase shaper (32).

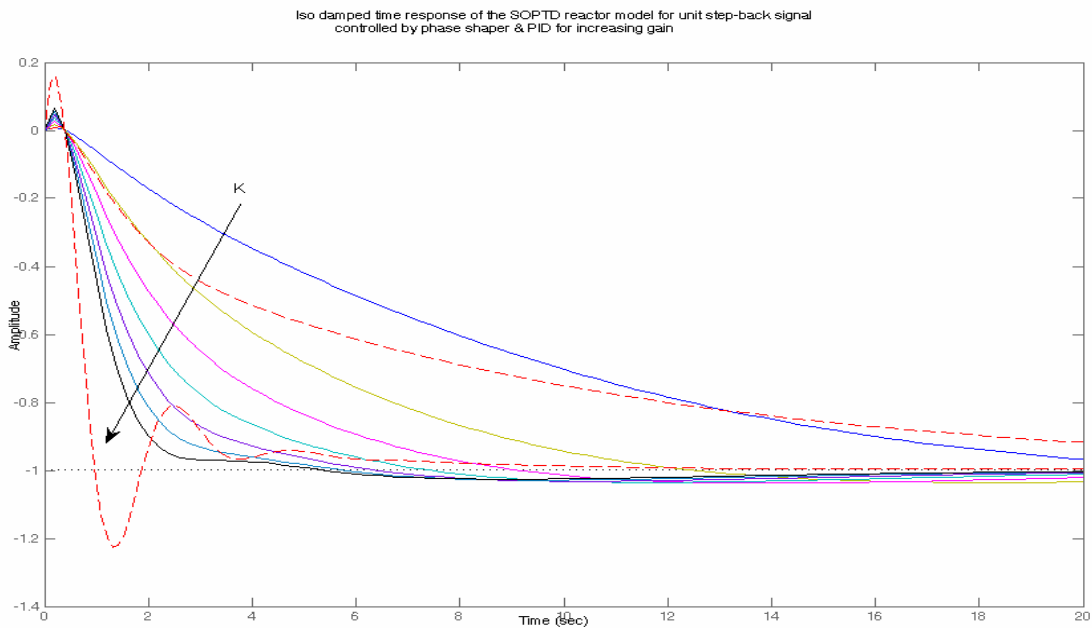

Fig. 5. Iso damped time response of the SOPTD reactor model (18) for unit step-back signal; dotted lines represent an active step-back response with only PID controller (31) and continuous lines represent the response of the active rod drop with the PID along with the phase shaper (32) for 7 times increase in gain.

From Fig. 5, it is again evident that using SOPDT approximation also, the controller and the designed phase shaper can be used to cater to an active step-back initiated at various reactor power levels, as in the case of FOPDT approximations.

## V. SIMULATION AND RESULTS FOR ACTIVE STEP-BACK IN THE REACTOR

In this section the simulation results for a control rod movement under a step-back condition, as represented in Fig.1 are presented. For simulation, the 30% rod drop case is simulated using FOPDT and SOPDT approximated plants representing the RRS and also compared to the systems represented by these plants augmented with respective PID controllers and phase shapers as reported in the previous sections.

With the proposed robust controller (phase shaper along with a PID), the dynamic responses of the FOPTD reactor model due to an active control rod drop up to 30% are shown in Fig. 6 and Fig. 7 for different starting reactor powers. This methodology is shown to perform much better over the present RRS in terms of handling sluggish response for truncated ramp input (Fig. 1).

As it is evident from Fig. 7, the present RRS in practice with proportional only controller shows steady-state error [2], due to the low proportional gain for low initial powers. As the gain increases, the offset is gradually minimized. It is well known that an integral element in the controller can force the steady-state offset to zero. In [4], [5] a PI/PID controller is proposed for removal of the offset by a driven mechanism. But the addition of a pole at the origin makes the overall system prone to have oscillatory response. A FO phase-shaper, on the other hand, makes it possible to have a faster response having good tracking behaviour and no undershoot (i.e. deadbeat response). In the proposed robust controller, there is an inherent integral term in the PID and the phase shaper itself, which removes the offset problem. This also ensures iso-damping property and hence parametric robustness of the system.

The response of the PID controller (31) along with a phase shaper (32) to the reactor dynamics due to an active control rod drop is shown in Fig. 6 and Fig. 7 for the SOPTD reactor models also, having different starting powers.

From Fig. 6 and Fig. 7, the response of the overall system considering FOPTD structure of the nonlinear reactor plant is proved to have a faster response, though both FOPTD and SOPTD modelling gives dead-beat response for same amount of increase in loop gain, if an active rod drop with the robust controller is done. So, for practical implementation in the test PHWR, the robust controller (phase shaper and PID) derived from a simple FOPTD structure is recommended.

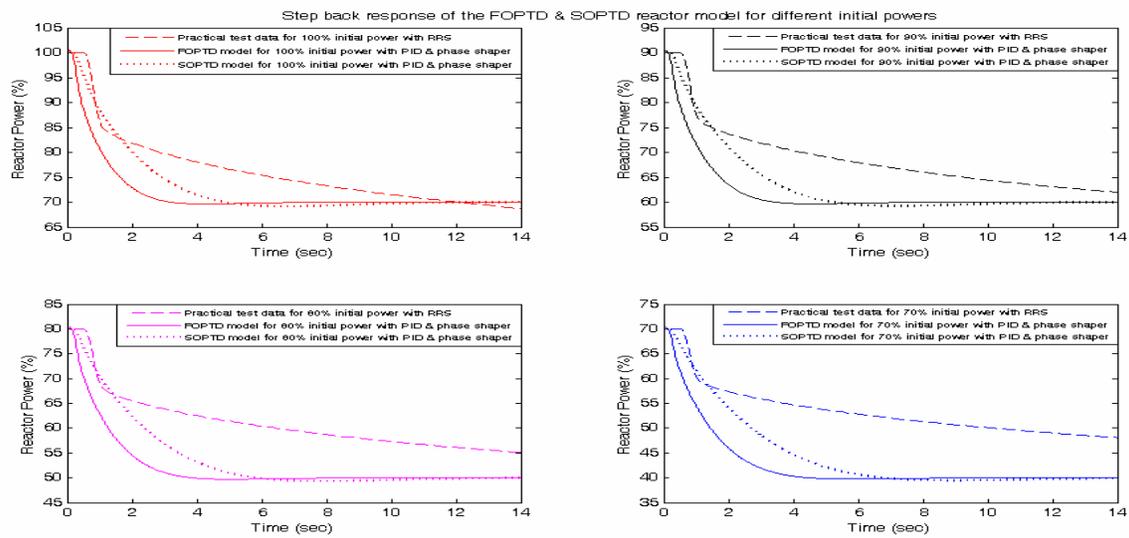

Fig. 6. 30% step-back response of the FOPTD and SOPTD reactor model operating at different initial powers.

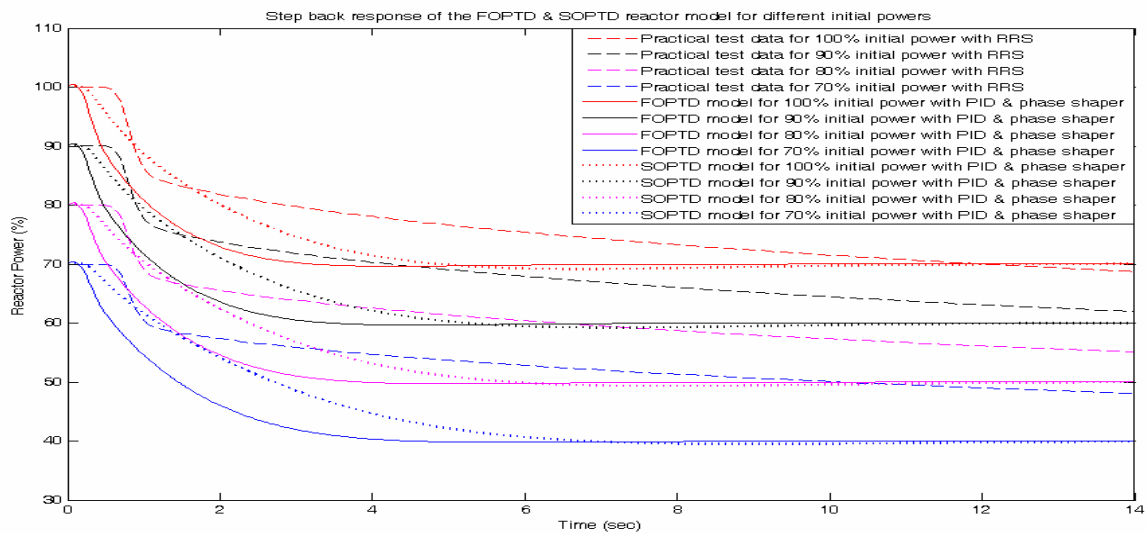

Fig. 7. Step-back response of the FOPTD and SOPTD model of the reactor at different initial power levels for 30% rod drop, continuous lines represent the response of the active rod drop with the PID along with the phase shaper for the two kind of modelling and dashed lines represent the dynamic response of passive rod drop with present RRS.

It is well known that the Sensitivity function $S(s)$ is an indication of the ability of the system to suppress load disturbances and achieve good set-point tracking while the complementary sensitivity function $T(s)$ indicates the robustness to measurement noise and other unmodelled system dynamics [33]. To have a good time response under these disturbed conditions, the $S(s)$ should have small values at lower frequencies while $T(s)$ should have small values at higher frequencies [23]. From Fig. 8, it is clear that for this system with a phase shaper, the value $S(j\omega)$ at low frequency ($\omega = 0.1\ rad/s$) is slightly increased and $T(j\omega)$ at high frequency ($\omega = 10\ rad/s$) is reduced considerably. The corresponding unit step-back and load disturbance responses are shown in Fig. 9. Also considerable reduction in initial controller output signal can be found in Fig. 10, due to the introduction of the phase shaper.

With the introduction of the phase shaper, significant improvements in robustness (in terms of gain variation), reduced controller outputs (less chance of actuator saturation and lesser actuator size) and complementary sensitivity at high frequency (better high frequency noise rejection) are obtained. But these three improvements are obtained at the cost of slight increase in sensitivity function at low frequency for the FOPTD reactor model and hence poor load disturbance response. This problem of inferior load disturbance response has been tackled in [23] with the inclusion of additional constraints limiting the values of $S(j\omega)$, while designing the fractional order controller itself. Also, in the load disturbance response (Fig. 9), the overshoot is slightly increased due to the lagging nature of the phase shaper and hence due to the reduction in phase margin. But as reported in (28), the optimization problem itself is designed to take care of the minimum desired phase margin ($\phi_{md}$) as a nonlinear constraint, which effectively controls the overshoot of the overall closed loop system.

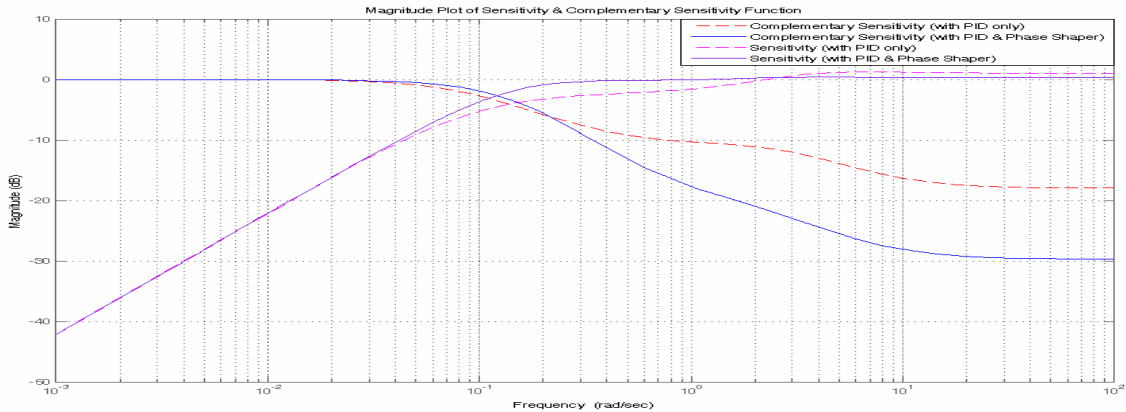

Fig. 8. Magnitude plot of sensitivity and complementary sensitivity due to phase shaping.

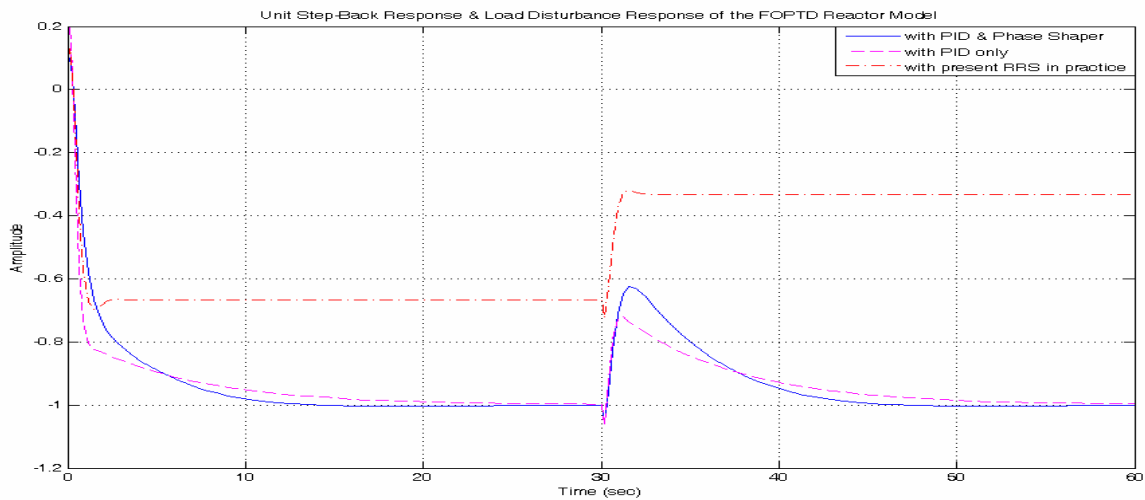

Fig. 9. Step-back response and load disturbance response of plant (14) with controller (31) with and without phase shaper (32).

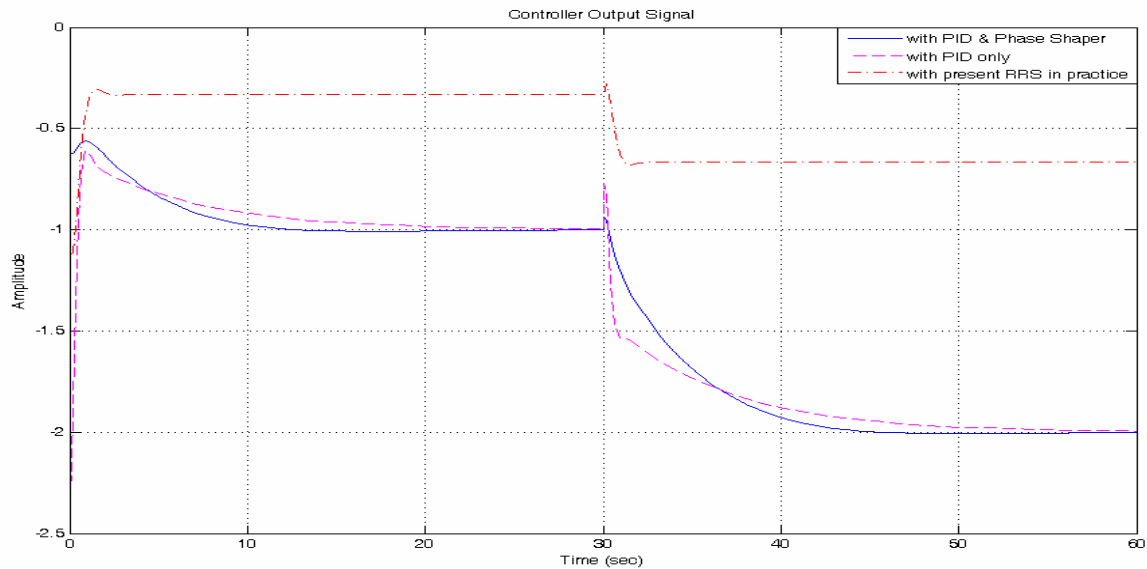

Fig. 10. Controller output signal for plant (14) with controller (31) with and without phase shaper (32).

For the present work, it is assumed that the PID controllers used for simulation are assumed to be tuned by the methodology described in [6] which guarantees a certain minimum phase margin. The constraint (28) specifying the minimum acceptable phase margin can be used to get desired closed-loop response. Tuning the PID control loops with a higher value of damping e.g. methods specifying gain and phase margins are likely to result in acceptable closed-loop response with iso-damping when the phase shaper is included in the loop.

While simulating the time and frequency response of the delay term of the plants, Pade's First Order approximation is used throughout the paper. For all the results presented, as mentioned in Section IV and V, $s^q$ is assumed to be represented by a first order Carlson's approximation which is a rational transfer function with a first order numerator and a first order denominator. The constrained optimization technique automatically establishes the parameters of the phase shaper $\{q,a\}$ that produces the maximum flat-phase frequency spread around the specified $\omega_{gc}$. With the $\{q,a\}$ obtained for each phase shaper, the results were repeated using higher order Carlson's representation and almost identical results were obtained. Higher order Carlson's representations of a fractional differ-integrator use rational transfer functions with higher order polynomials as integer order approximations of a FO differ-integrator with a greater accuracy. This proves the adequacy of a first order Carlson approximated FO differ-integrator and shows that the methodology presented in this paper can be used to design a practically realizable phase shaper which can be used in conjunction with the RRS of the presents PHWRs.

## VI. Conclusion

In this paper, an active step-back mechanism has been proposed with a robust controller comprising of a FO phase shaper and a LQR tuned PID. The methodology, put forward in this paper, is proved to be far better than the RRS in practice for the Indian PHWRs as far as parametric robustness against ageing, nonlinearity and large undershoot for high system gains, are of major concern. The methodology proposed in this paper, can be used for enhancing the parametric robustness of the RRS control loop with gain variations, while achieving considerable fast response with excellent tracking behaviour. The iso-damped nature of the response allows design of extremely fast step-back mechanism, keeping the overshoot constant, provided the actuator constraints can be met. Further, the phase shaper proposed is of a low-order and practically realizable. Thus, in practical terms the methodology allows design of a simple hardware element which can be used with PID controller, tuned by any standard method.

Modelling of the nonlinear reactor as a fractional order plus time delay system and an active step-back mechanism with a $PI^\lambda D^\mu$ controller to stabilize the time-delay as well as the complicated fractional order dynamics are a bit challenging and is left as the scope for future work.

## Acknowledgment

The work presented in this paper has been supported by the Board of Research in Nuclear sciences (BRNS) of the Department of Atomic Energy, India. Sanction no 2006/34/34-BRNS dated March 2007.